\documentclass[11pt,a4paper]{article}

\usepackage[utf8]{inputenc}
\usepackage[T1]{fontenc}
\usepackage[english]{babel}
\usepackage{amsmath,amssymb,amsthm}
\usepackage{mathtools}
\usepackage{geometry}
\usepackage{booktabs}
\usepackage{array}
\usepackage{enumitem}
\usepackage{lmodern}
\usepackage{hyperref}
\hypersetup{
  colorlinks=true,
  linkcolor=blue,
  citecolor=blue,
  urlcolor=blue,
  pdftitle={Perimeter-Regularized Critical Shapes in the Class O_C},
  pdfauthor={M. Barkatou, M. El Morsalani}
}
\usepackage{microtype}

\theoremstyle{plain}
\newtheorem{theorem}{Theorem}[section]
\newtheorem{proposition}[theorem]{Proposition}
\newtheorem{lemma}[theorem]{Lemma}
\newtheorem{corollary}[theorem]{Corollary}
\newtheorem{conjecture}[theorem]{Conjecture}

\theoremstyle{definition}
\newtheorem{definition}[theorem]{Definition}
\newtheorem{assumption}[theorem]{Assumption}
\newtheorem{remark}[theorem]{Remark}

\newcommand{\R}{\mathbb{R}}

\newcommand{\Sph}{\mathbb{S}}

\newcommand{\DR}{D\mathcal{R}}
\newcommand{\Gtil}{\widetilde{G}}
\newcommand{\norm}[1]{\left\lVert #1 \right\rVert}
\newcommand{\abs}[1]{\left\lvert #1 \right\rvert}

\newcommand{\dd}{\,\mathrm{d}}
\newcommand{\eps}{\varepsilon}
\newcommand{\supp}{\operatorname{supp}}
\newcommand{\intC}{\operatorname{int}C}
\newcommand{\diverg}{\operatorname{div}}
\newcommand{\Ret}{\Pi^{\mathrm{ret}}_0}

\title{%
  \textbf{Perimeter-Regularized Critical Shapes in the Class $\mathcal{O}_C$:}\\[2pt]
  \textbf{Core-Normal Transfer and First-Order Return-Map Response}%
}
\author{M.~Barkatou\thanks{Corresponding author.}\quad M.~El Morsalani}
\date{Research draft v0.5 --- November 2026}

\begin{document}
\maketitle

% ============================================================
% ABSTRACT
% ============================================================
\begin{abstract}
Let $C\subset\R^N$ be a compact strictly convex core and let $\mathcal{O}_C$ denote the class of domains whose outer boundary is a global outward normal graph over $\partial C$ and whose inward boundary normals return to the interior of $C$. We consider volume-constrained critical domains of the Dirichlet-perimeter functional
\[
\mathcal{E}_{\mu}(\Omega) = -\frac{1}{2}\int_{\Omega}\abs{\nabla u_{\Omega}}^{2}\dd x + \mu P(\Omega),
\qquad -\Delta u_{\Omega} = f,\quad u_{\Omega}|_{\partial\Omega}=0,
\]
in the large-perimeter regime $\mu\to\infty$. The local existence mechanism near translated balls is classical in shape optimization, going back in particular to Hayouni--Novruzi and Novruzi. Our purpose is different: to connect such a critical-shape branch to the normal-return geometry of $\mathcal{O}_C$ and to quantify the induced deformation of its return map.

We first give a self-contained adaptation of the translated-ball Lyapunov--Schmidt construction in a form that extracts the first physical shape coefficient. If $x_0$ is a nondegenerate critical point of the reduced translated-ball energy, the branch has genuine shape scale $\mu^{-1}$ and
\[
\psi_{\mu} = \mu^{-1}\psi_{1} + o(\mu^{-1}),
\qquad J_{R}\psi_{1} = \tfrac{1}{2} P_{X}(g_{0}^{2}).
\]
The new geometric part then begins. We prove that $\mathcal{O}_C$ is open under sufficiently small $C^1$ perturbations of the thickness, establish a quantitative coordinate-transfer theorem from the selected sphere to a unique core-normal thickness, and prove that the forward return-map operator
\[
\mathcal{R}:d\longmapsto F_{d}
\]
is Fréchet $C^1$ between Hölder spaces. Its derivative is obtained explicitly from the first-hit equation. Combining these results yields
\[
d_{\mu} = d_{0} + \mu^{-1}\eta_{1} + o(\mu^{-1}),
\qquad
F_{\mu} = F_{0} + \mu^{-1}\cal{R}_d[d_{0}]\eta_{1} + o(\mu^{-1}).
\]
Thus the paper supplies a shape-to-dynamics sensitivity theorem for the class $\mathcal{O}_C$: a perimeter-regularized PDE-constrained critical deformation produces an explicitly controlled first-order deformation of the associated normal-return dynamics. A perturbative nonradial example shows that the response coefficient $\eta_1$ is generically nonzero. No claim of global minimality is made.
\end{abstract}



% ============================================================
% 1. INTRODUCTION
% ============================================================
\section{Introduction}

\subsection{Problem, positioning, and guiding obstruction}
Let $C\subset\R^N$, $N\geq 2$, be a compact strictly convex body. The class $\mathcal{O}_{C}$, originating in the geometric normal-property framework of Barkatou, consists here of domains surrounding $C$ whose outer boundary is a global normal graph over $\partial C$ and whose inward outer-boundary normals return to $\intC$; see \cite{Barkatou2002,OC1,OC3}. The geometry carries a return map
\[
F_{d}:\partial C\longrightarrow\partial C,
\]
obtained by moving first along the outward core normal to the outer boundary and then back to the core along the inward normal of that outer boundary.

We study the interaction of this return geometry with the volume-constrained Dirichlet-perimeter functional
\begin{equation}\label{eq:energy}
\mathcal{E}_{\mu}(\Omega) = -\frac{1}{2}\int_{\Omega}\abs{\nabla u_{\Omega}}^{2}\dd x + \mu P(\Omega),
\qquad \mu>0,
\end{equation}
where $u_{\Omega}$ solves
\begin{equation}\label{eq:state}
-\Delta u_{\Omega} = f\quad\text{in }\Omega,
\qquad u_{\Omega}=0\quad\text{on }\partial\Omega,
\end{equation}
and $\nu_{\Omega}$ denotes the outward unit normal. We impose
\begin{equation}\label{eq:volume}
\abs{\Omega} = V_0.
\end{equation}
The source $f$ is fixed in physical space.

The large-$\mu$ regime is dominated by perimeter. A necessary-limit argument shows that any connected smooth $C^2$ limit of critical domains with uniformly bounded boundary flux must have constant mean curvature and is therefore spherical. Thus the correct reference object is a selected ball of volume $V_{0}$, not an arbitrary nonspherical core. The convex core remains essential in a different role: it provides the global normal coordinates and the return geometry through which the critical branch is represented and observed.

\subsection{What is inherited from shape optimization and what is new here}
The translated-ball existence mechanism itself is not new. Hayouni and Novruzi \cite{HayouniNovruzi} used the implicit-function theorem to construct regular critical points for a two-dimensional Dirichlet-perimeter free-boundary problem for large perimeter weight. Novruzi \cite{Novruzi} subsequently developed a $C^{2,\alpha}$ existence theory for a broader class of shape functionals, explicitly separating translations from genuine shape deformations and using nondegeneracy of the reduced translated-ball energy. Related free-boundary work by Barkatou, Seck and Ly \cite{BarkatouSeckLy} treated an interior electromagnetic-casting problem through shape optimization and curvature boundary conditions.

Accordingly, the Lyapunov--Schmidt construction in Sections~\ref{sec:prep}--\ref{sec:LS} is included as a self-contained preparatory layer adapted to the present notation and to the first-order coefficient required later. The principal new results begin when this classical branch is connected to $\mathcal{O}_{C}$:
\begin{itemize}[leftmargin=2em]
\item[(N1)] a $C^{1}$-openness theorem for the geometric normal-return class $\mathcal{O}_{C}$ (Theorem~\ref{thm:open});
\item[(N2)] a quantitative coordinate transfer from a translated spherical graph to a unique core-normal thickness, including the joint contribution of center drift and higher spherical modes (Theorem~\ref{thm:transfer});
\item[(N3)] Fréchet $C^{1}$ differentiability of the forward return operator $\mathcal{R}:d\mapsto F_{d}$, with an explicit derivative obtained from the first-hit equation (Theorem~\ref{thm:DR}).
\end{itemize}
The paper therefore does not claim a new general existence theory for perimeter-regularized critical shapes. Its novelty is the shape-to-return-dynamics sensitivity theory in the class $\mathcal{O}_C$.

\subsection{Relation with the preceding $\mathcal{O}_C$ papers}
The present paper is the next step in a sequence with a progressively changing object of study. The first return-map paper \cite{OC1} established the exact normal geometry of $d\mapsto F_{d}$, the fixed-point/critical-point relation, and the spatial linearization $D_{c}F_{d}$ at fixed points. The inverse paper \cite{OC2} studies the opposite direction, reconstructing thickness information from a given return map through an exact nonlinear transport system and its linearized compatibility structure. The elliptic superlevel-set paper \cite{OC3} keeps the PDE domain fixed and studies one distinguished deformation path
\[
d_{t} = d_{0} - t b_{0} + o_{C^{1}}(t),
\qquad
F_{t} = F_{0} + tA_{0} + o(t),
\]
generated by small Dirichlet superlevels.

The present problem is complementary to all three. The domain itself is varied by an Euler--Lagrange condition, and the main forward object is now the operator
\[
\mathcal{R}:d\longmapsto F_{d}.
\]
Theorem~\ref{thm:DR} differentiates this operator with respect to an arbitrary thickness perturbation. Hence the earlier PDE coefficient $A_{0}$ becomes one special direction of a general Fréchet shape derivative, while the inverse-paper compatibility system becomes a natural target for a future range characterization of $\cal{R}_d[d_0]$.

\subsection{Comparison with recent work}
Recent shape-optimization literature emphasizes several complementary questions. Buttazzo--Maiale--Mazzoleni--Tortone--Velichkov \cite{Buttazzo} develop free-boundary regularity for broad PDE-dependent integral shape functionals. Celentano--Nitsch--Trombetti \cite{Celentano} use overdetermined elliptic systems to characterize critical shapes under volume constraints. Caubet--Dambrine--Darde \cite{Caubet} show that perimeter penalization can improve stability in a Dirichlet inverse-obstacle shape problem. Plotnikov--Sokolowski \cite{Plotnikov} study geometric gradient-flow frameworks for regularized inverse and free-boundary shape problems, and Amato--Nitsch--Trombetti--Villone \cite{Amato} give recent local stability results for functionals combining Dirichlet spectral/torsional quantities with perimeter constraints.

These works concern existence, regularity, stability, symmetry, or optimization dynamics of the shape itself. The question addressed here is different: once a smooth critical-shape branch is available, how does an additional geometric observable---the normal return map attached to a convex core---change with the shape? The new object $\DR[d_0]$ is therefore a response operator from boundary deformation to discrete return dynamics rather than a shape Hessian or a gradient-flow generator.

Two adjacent theories should be distinguished from this claim. In classical Poincaré/first-impact theory, transversality of a flow with a section yields differentiability of the return time and hence of the first-return map by an implicit-function argument; see, for example, Immler--Traut \cite{ImmlerTraut}. Thus the projection that removes the component generated by variation of the hitting time is not, by itself, new. In shape calculus, Feppon--Allaire--Dapogny \cite{Feppon} study shape derivatives of geometric constraints involving signed distance and integration along normal rays. Their normal-ray computations are closely related to our equations \eqref{eq:normals}--\eqref{eq:delta-n}, but they do not address return maps.

\subsection{Theorem hierarchy and main result}
To make the novelty hierarchy explicit, we distinguish one preparatory theorem from the principal geometric results.

\begin{table}[ht]
\centering
\small
\renewcommand{\arraystretch}{1.25}
\begin{tabular}{@{}p{3.6cm}p{4.6cm}p{5.2cm}@{}}
\toprule
\textbf{Work} & \textbf{Main issue} & \textbf{Relation to the present paper}\\
\midrule
Hayouni--Novruzi (2002); Novruzi (2004) & Large-perimeter critical/optimal domains; translation versus pure deformation & Classical existence mechanism used as preparatory background; no $\mathcal{O}_C$ return response\\
Barkatou--Seck--Ly (2006) & Curvature free-boundary condition and electromagnetic casting & Motivates the PDE-perimeter model and applications; no normal-return sensitivity\\
Buttazzo et al.\ (2024) & Regularity of PDE-dependent optimal sets & Global/free-boundary regularity, complementary to our local smooth branch\\
Caubet--Dambrine--Darde (2026) & Perimeter regularization in inverse obstacle optimization & Shows current relevance of perimeter stabilization; different observable and inverse problem\\
Earlier $\mathcal{O}_C$ papers & Geometry, inverse reconstruction, and one PDE-generated deformation of $F_d$ & Supply the exact return geometry; the present work differentiates the full forward operator $d\mapsto F_d$\\
\textbf{Present paper} & Core-normal transfer and first-order return-map response & Fréchet response operator and combined critical-shape-to-dynamics theorem\\
\bottomrule
\end{tabular}
\caption{Position of the present results relative to closely related work.}
\label{tab:position}
\end{table}

\begin{enumerate}[leftmargin=2em]
\item[(1)] \emph{Preparatory selected-ball branch.} We adapt the classical translated-ball implicit-function/Lyapunov--Schmidt framework and extract the physical $\mu^{-1}$ coefficient
\[
J_{R}\psi_{1} = \tfrac{1}{2} P_{X}(g_{0}^{2}).
\]
\item[(2)] \emph{$\mathcal{O}_C$ openness and transfer.} Small $C^1$ thickness perturbations preserve $\mathcal{O}_C$, and a nearby translated spherical graph has a unique core-normal representation.
\item[(3)] \emph{Principal response theorem.} The forward return operator is Fréchet $C^1$, with an explicit derivative $\DR[d_0]$.
\item[(4)] \emph{Combined shape-to-dynamics theorem.} Along the critical branch,
\[
d_{\mu} = d_{0} + \mu^{-1}\eta_{1} + o_{C^{2,\alpha}}(\mu^{-1}),
\qquad
F_{\mu} = F_{0} + \mu^{-1}\DR[d_{0}]\eta_{1} + o_{C^{1,\alpha}}(\mu^{-1}).
\]
\end{enumerate}
No theorem below uses numerical evidence. Exact geometric identities are established before asymptotic expansions; every asymptotic statement is labelled by its topology and small parameter.

\subsection{Potential applications}
The applications are secondary to the geometric analysis, but the theorem chain has several natural interpretations. First, Dirichlet-perimeter free-boundary conditions arise historically in electromagnetic casting and shaping, where curvature regularization competes with a PDE-generated boundary pressure; see \cite{HayouniNovruzi,BarkatouSeckLy}. Second, perimeter penalties remain active in contemporary inverse-shape problems because they can improve stability and numerical regularity \cite{Caubet}. Third, the new operator $\DR[d_0]$ provides a sensitivity map from boundary design variables to a geometric observable. This can be relevant whenever an interface is probed or controlled through normal incidence/return information, for example in inverse geometry, geometric metrology, or PDE-constrained interface design. We do not claim an application-specific sensing model here; the contribution is the mathematical response operator that such models could use.

\subsection{Methodological conventions and reproducibility}
Several conventions are fixed globally to avoid ambiguities that are particularly damaging in normal-return problems. We use
\[
S = -D\nu|_{TM},
\]
so that a convex core has $S\leq 0$ and $I-dS\geq I$ for every outward thickness $d\geq 0$. Outward and inward normals are denoted by $N_{d}$ and $n_{d}=-N_{d}$, respectively. Exact identities, first-order expansions, and heuristic interpretations are kept separate. The analytical benchmark in Section~\ref{sec:bench} is exact. A computational protocol is stated there for a nonradial test, but no numerical result is used as evidence for any theorem.

\subsection{Scope and nonclaims}
The results are local. We construct constrained critical domains; we do not claim global minimality. The existence of a selected nondegenerate center and of a containing reference sphere of the prescribed volume is assumed (see Remark~\ref{rem:H3}). We do not claim that $\operatorname{Ran}\DR[d_0]$ equals the full compatibility space of the linearized inverse problem; a precise conjecture is stated in Section~\ref{sec:open}. Finally, the preparatory large-perimeter branch is presented as an adaptation of established shape-optimization machinery, not as a new general existence theorem.

% ============================================================
% 2. GEOMETRIC AND VARIATIONAL SETUP
% ============================================================
\section{Geometric and variational setup}\label{sec:setup}

\subsection{The class $\mathcal{O}_C$}
Let
\[
M\coloneqq\partial C\in C^{4,\alpha},
\qquad 0<\alpha<1,
\]
and let $\nu:M\to\Sph^{N-1}$ be the outward unit normal. We adopt the convention
\begin{equation}\label{eq:S}
S = -D\nu|_{TM}.
\end{equation}
For a convex hypersurface the eigenvalues of $S$ are $-\kappa_{1},\dots,-\kappa_{N-1}$, hence $S\leq 0$.

For $d:M\to(0,\infty)$ define
\begin{equation}\label{eq:Psi-d}
\Psi_d(c) = c + d(c)\nu(c),
\end{equation}
and
\begin{equation}\label{eq:Omega-d}
\Omega[d] = C\cup\{c+r\nu(c):c\in M,\ 0<r<d(c)\}.
\end{equation}

\begin{definition}[The class $\mathcal{O}_C$]\label{def:OC}
A domain $\Omega[d]$ with $d\in C^{2,\alpha}(M)$ and $d>0$ belongs to $\mathcal{O}_{C}$ if every inward normal ray issued from $\partial\Omega[d]$ intersects $\intC$.
\end{definition}

The global normal-graph property is already built into \eqref{eq:Omega-d}. Set
\begin{equation}\label{eq:w-chi}
w_{d} = L_{d}^{-1}\nabla_{M}d,
\qquad
\chi_{d} = \sqrt{1+\abs{w_{d}}^{2}}.
\end{equation}
The exact outward and inward unit normals are
\begin{equation}\label{eq:normals}
N_{d} = \frac{\nu-w_{d}}{\chi_{d}},
\qquad
n_{d} = -N_{d} = \frac{-\nu+w_{d}}{\chi_{d}}.
\end{equation}
These formulas are exact and are the starting point for all later differentiations; compare \cite{OC3}.

\subsection{Energy and Euler--Lagrange equation}
Fix
\begin{equation}\label{eq:f}
f\in C_{c}^{4,\alpha}(\R^{N}),
\qquad f\geq 0,
\qquad f\neq 0.
\end{equation}
For each smooth domain containing $\supp f$, let $u_{\Omega}$ solve \eqref{eq:state}.

For a smooth deformation with normal velocity $V_{N}$, the Hadamard formula (see, e.g., \cite{DelfourZolesio,HenrotPierre,SokolowskiZolesio}) gives
\begin{equation}\label{eq:hadamard}
D\!\left(-\frac{1}{2}\int_{\Omega}\abs{\nabla u_{\Omega}}^{2}\right)[V_{N}]
= -\frac{1}{2}\int_{\partial\Omega}(\partial_{N}u_{\Omega})^{2}V_{N}\dd S.
\end{equation}
With the convention $H=\diverg_{\partial\Omega}N$, the perimeter and volume derivatives are
\[
DP(\Omega)[V_N] = \int_{\partial\Omega}HV_N\dd S,
\qquad
D\abs{\Omega}[V_N] = \int_{\partial\Omega}V_N\dd S.
\]
Thus a smooth volume-constrained critical domain satisfies
\begin{equation}\label{eq:EL}
-\tfrac{1}{2}g^{2} + \mu H + \Lambda_{\mu} = 0,
\qquad
g\coloneqq-\partial_{N}u_{\Omega}.
\end{equation}
Set
\begin{equation}\label{eq:eps}
\eps = \mu^{-1},
\qquad
\lambda_{\eps} = \Lambda_{\mu}/\mu.
\end{equation}
Then
\begin{equation}\label{eq:EL-scaled}
H + \lambda_{\eps} - \tfrac{\eps}{2}g^{2} = 0.
\end{equation}

\subsection{Rigidity of a smooth large-$\mu$ limit}
\begin{proposition}[CMC obstruction]\label{prop:CMC}
Let $\Omega_{\mu_{j}}$ be smooth volume-constrained critical domains for \eqref{eq:energy}, with $\mu_{j}\to\infty$. Assume
\[
\partial\Omega_{\mu_{j}}\to\Sigma_{0}\quad\text{in }C^{2}
\]
and
\[
\sup_{j}\norm{\partial_{N_{j}}u_{\mu_{j}}}_{L^{\infty}(\partial\Omega_{\mu_{j}})}<\infty.
\]
Then $\Sigma_{0}$ has constant mean curvature.
\end{proposition}

\begin{proof}
Subtract the surface average from \eqref{eq:EL}:
\[
\mu_{j}(H_{j}-\overline{H}_{j}) = \tfrac{1}{2}(g_{j}^{2}-\overline{g_{j}^{2}}).
\]
The flux bound implies
\[
\norm{H_{j}-\overline{H}_{j}}_{L^{\infty}}\leq\frac{C}{\mu_{j}}\to 0.
\]
The $C^2$ convergence gives uniform convergence of the mean curvatures, so $H_{\Sigma_0}$ is constant.
\end{proof}

\begin{remark}\label{rem:alexandrov}
In the compact connected embedded setting relevant here, a closed constant-mean-curvature hypersurface is spherical by Alexandrov's theorem \cite{Alexandrov}. Proposition~\ref{prop:CMC} is only a necessary-limit statement; it does not by itself provide compactness or existence of a large-$\mu$ branch.
\end{remark}

% ============================================================
% 3. SELECTED SPHERICAL REFERENCE STATES
% ============================================================
\section{Selected spherical reference states: preparatory reduction}\label{sec:prep}

Fix $R>0$ by
\begin{equation}\label{eq:R}
\abs{B_R} = V_0.
\end{equation}
Let $U\subset\R^{N}$ be a neighborhood of a point $x_0$ such that
\begin{equation}\label{eq:supp}
\supp f\Subset B_{R-\delta}(x)\qquad\text{for every }x\in U
\end{equation}
for some $\delta>0$.

For $x\in U$, let $u_{x}$ solve
\[
-\Delta u_{x} = f\quad\text{in }B_{R}(x),
\qquad u_{x}=0\quad\text{on }\partial B_{R}(x).
\]
Define the translated-ball reduced energy
\begin{equation}\label{eq:h}
h(x) = -\frac{1}{2}\int_{B_R(x)}\abs{\nabla u_x}^{2}\dd y.
\end{equation}

\begin{lemma}[Gradient of the reduced energy]\label{lem:grad-h}
Let
\[
g_{x}(\omega) = -\partial_{N}u_{x}(x+R\omega).
\]
Then
\begin{equation}\label{eq:grad-h}
\nabla h(x) = -\frac{R^{N-1}}{2}\int_{\Sph^{N-1}}g_{x}(\omega)^{2}\omega\dd\omega.
\end{equation}
\end{lemma}

\begin{proof}
Translate the ball in the direction $a\in\R^{N}$. The normal velocity is $a\cdot\omega$. The Hadamard formula \eqref{eq:hadamard} gives
\[
Dh(x)[a] = -\frac{1}{2}\int_{\partial B_{R}(x)}g_{x}^{2}(a\cdot N)\dd S,
\]
which is equivalent to \eqref{eq:grad-h}.
\end{proof}

\begin{assumption}[Selected nondegenerate center]\label{ass:center}
There exists $x_0\in U$ such that
\begin{equation}\label{eq:nondeg}
\nabla h(x_0) = 0,
\qquad
\det D^{2}h(x_0)\neq 0.
\end{equation}
\end{assumption}

The nondegeneracy is the finite-dimensional hypothesis needed to resolve the translation kernel of the spherical Jacobi operator.

% ============================================================
% 4. EXACT-VOLUME CHART AND SMOOTH PDE DEPENDENCE
% ============================================================
\section{Exact-volume chart and smooth PDE dependence}\label{sec:chart}

Let $\sigma_{N-1}=\abs{\Sph^{N-1}}$. Define
\begin{equation}\label{eq:X}
X^{k,\alpha} = \left\{\psi\in C^{k,\alpha}(\Sph^{N-1}):\int_{\Sph^{N-1}}\psi\dd\omega=0,\ \int_{\Sph^{N-1}}\psi\,\omega\dd\omega=0\right\}.
\end{equation}

\subsection{Exact volume preservation}
For $x\in U$ and $\psi\in X^{3,\alpha}$, consider
\begin{equation}\label{eq:Sigma-psi}
\Sigma_{x,\psi} = \{x+[R+a(\psi)+\psi(\omega)]\omega:\omega\in\Sph^{N-1}\}.
\end{equation}

\begin{lemma}[Exact volume chart]\label{lem:volume}
There is a neighborhood of $0$ in $X^{3,\alpha}$ and a unique smooth scalar function $a(\psi)$ such that the domain bounded by \eqref{eq:Sigma-psi} has volume $V_{0}$. Moreover
\begin{equation}\label{eq:a0}
a(0)=0,
\qquad Da(0)=0,
\end{equation}
and
\begin{equation}\label{eq:a2}
a(\psi) = -\frac{N-1}{2R\sigma_{N-1}}\int_{\Sph^{N-1}}\psi^{2}\dd\omega + O(\norm{\psi}_{C^{3,\alpha}}^{3}).
\end{equation}
\end{lemma}

\begin{proof}
The enclosed volume is
\[
V(a,\psi) = \frac{1}{N}\int_{\Sph^{N-1}}(R+a+\psi)^{N}\dd\omega.
\]
At $(0,0)$,
\[
\partial_{a}V(0,0) = R^{N-1}\sigma_{N-1}>0.
\]
The implicit-function theorem gives $a(\psi)$. Since $\psi$ has zero mean,
\[
D_{\psi}V(0,0)[\phi] = R^{N-1}\int_{\Sph^{N-1}}\phi = 0,
\]
so $Da(0)=0$. Differentiating twice gives
\[
D^{2}a(0)[\phi,\eta] = -\frac{N-1}{R\sigma_{N-1}}\int_{\Sph^{N-1}}\phi\eta,
\]
which yields \eqref{eq:a2}.
\end{proof}

\subsection{Smooth dependence of the state and flux}
\begin{lemma}[Smooth state and flux map]\label{lem:smooth}
Under \eqref{eq:f}--\eqref{eq:supp}, after shrinking $U$ and the neighborhood of $0$ in $X^{3,\alpha}$ if necessary, the solution $u_{x,\psi}$ of
\[
-\Delta u_{x,\psi} = f\quad\text{in }\Omega_{x,\psi},
\qquad u_{x,\psi}=0\quad\text{on }\partial\Omega_{x,\psi},
\]
depends $C^{2}$-smoothly on $(x,\psi)$ after pullback to $B_{R}$. The pulled-back boundary flux
\begin{equation}\label{eq:flux}
g(x,\psi)(\omega) = -\partial_{N}u_{x,\psi}(x+[R+a(\psi)+\psi(\omega)]\omega)
\end{equation}
is a $C^{2}$ map into $C^{1,\alpha}(\Sph^{N-1})$.
\end{lemma}

\begin{proof}
Choose a smooth radial cutoff $\chi$ equal to $0$ on $[0,R-2\delta]$ and $1$ near $R$. Define $T_{x,\psi}(y)=x+y$ on a fixed neighborhood of the origin and, writing $y=r\omega$ outside that neighborhood,
\[
T_{x,\psi}(r\omega) = x+[r+\chi(r)(a(\psi)+\psi(\omega))]\omega.
\]
Because the two formulas agree where $\chi=0$, this gives a globally $C^{3,\alpha}$ map on $B_{R}$. For $\psi$ small it is a $C^{3,\alpha}$ diffeomorphism from $B_{R}$ onto $\Omega_{x,\psi}$.

Put
\[
X_{D} = C_{0}^{2,\alpha}(B_{R}),
\qquad Y = C^{0,\alpha}(B_{R}),
\]
and define
\[
\mathcal{F}(x,\psi,v) = -\diverg(A_{\psi}\nabla v) - J_{\psi}f(T_{x,\psi}),
\]
with
\[
A_{\psi} = J_{\psi}DT_{x,\psi}^{-1}DT_{x,\psi}^{-T},
\qquad J_{\psi} = \det DT_{x,\psi}.
\]
Then
\[
\mathcal{F}:U\times X^{3,\alpha}\times X_{D}\longrightarrow Y
\]
is $C^{2}$: determinant and matrix inversion are smooth on the uniformly invertible coefficient set, the composition $(x,\psi)\mapsto f\circ T_{x,\psi}$ is $C^{2}$ because $f\in C_{c}^{4,\alpha}$, and multiplication/divergence have the required Hölder mapping properties. Moreover,
\[
D_{v}\mathcal{F}(x,\psi,v)z = -\diverg(A_{\psi}\nabla z)
\]
is a uniformly elliptic Dirichlet isomorphism $X_{D}\to Y$ by Schauder theory. The Banach-space implicit-function theorem therefore yields a $C^{2}$ map $(x,\psi)\mapsto v_{x,\psi}\in C_{0}^{2,\alpha}(B_{R})$.

Finally,
\[
\nabla u_{x,\psi} = DT_{x,\psi}^{-T}\nabla v_{x,\psi},
\]
and the outer normal of $\Sigma_{x,\psi}$ depends $C^{2}$-smoothly on $\psi$. The trace map $C^{2,\alpha}(B_{R})\to C^{1,\alpha}(\Sph^{N-1})$ is bounded. Hence the pulled-back flux \eqref{eq:flux} is $C^{2}$ into $C^{1,\alpha}(\Sph^{N-1})$.
\end{proof}

% ============================================================
% 5. ADAPTED LYAPUNOV-SCHMIDT CONSTRUCTION
% ============================================================
\section{Adapted Lyapunov--Schmidt construction}\label{sec:LS}

Let $H(\psi)$ denote the mean curvature of $\Sigma_{x,\psi}$ pulled back to $\Sph^{N-1}$. It is independent of $x$ by translation invariance. Define
\begin{equation}\label{eq:M}
\mathcal{M}(x,\psi,\eps) = H(\psi) - \frac{\eps}{2}g(x,\psi)^{2}.
\end{equation}
Equation \eqref{eq:EL-scaled} holds for some scalar $\lambda_{\eps}$ exactly when the nonconstant part of $\mathcal{M}$ vanishes. Let
\[
T = \operatorname{span}\{\omega_{1},\dots,\omega_{N}\}
\]
be the translation space, let $P_{T}$ be the $L^{2}(\Sph^{N-1})$ projection onto $T$, and let $P_{X}$ be the projection onto the orthogonal complement of constants and $T$.

\begin{lemma}[Smooth mean-curvature map]\label{lem:H}
For $\psi$ in a sufficiently small neighborhood of $0$ in $X^{3,\alpha}$, the pulled-back mean-curvature map
\[
\psi\longmapsto H(\psi)
\]
is $C^{\infty}$ from $X^{3,\alpha}$ to $C^{1,\alpha}(\Sph^{N-1})$. Its differential at the sphere is
\[
DH(0)[\phi] = J_{R}\phi,
\qquad
J_{R} = -\Delta_{\partial B_{R}} - \frac{N-1}{R^{2}}.
\]
\end{lemma}

\begin{proof}
The metric, second fundamental form, inverse metric and unit normal of a radial graph depend smoothly on $\psi$, $\nabla\psi$ and $\nabla^{2}\psi$ as long as $R+a(\psi)+\psi$ stays positive. Since $a:X^{3,\alpha}\to\R$ is smooth, the resulting mean curvature is a smooth quasilinear second-order operator from $C^{3,\alpha}$ to $C^{1,\alpha}$. The standard first-variation formula for mean curvature at a sphere gives the displayed Jacobi operator with the convention $H=\diverg N$.
\end{proof}

\subsection{Inversion of the higher spherical modes}
At $\psi=0$,
\[
H(0) = \frac{N-1}{R}.
\]
The linearized mean-curvature operator is
\begin{equation}\label{eq:JR}
J_{R} = -\Delta_{\partial B_{R}} - \frac{N-1}{R^{2}} = -\frac{1}{R^{2}}\bigl(\Delta_{\Sph^{N-1}}+N-1\bigr).
\end{equation}
On degree-$\ell$ spherical harmonics,
\begin{equation}\label{eq:eigen}
J_{R}Y_{\ell} = \frac{\ell(\ell+N-2)-(N-1)}{R^{2}}Y_{\ell}.
\end{equation}
Thus constants form the volume mode, degree one is the translation kernel, and $J_{R}$ is positive and invertible on all degrees $\ell\geq 2$.

\begin{proposition}[Higher-mode equation]\label{prop:higher}
There is a neighborhood of $(x_0,0)$ and a unique $C^2$ map
\[
\Psi:U\times(-\eps_0,\eps_0)\to X^{3,\alpha}
\]
such that
\begin{equation}\label{eq:Psi-eq}
P_{X}\mathcal{M}(x,\Psi(x,\eps),\eps) = 0,
\qquad \Psi(x,0)=0.
\end{equation}
Moreover
\begin{equation}\label{eq:Psi1}
J_{R}\partial_{\eps}\Psi(x,0) = \tfrac{1}{2}P_{X}(g_{x}^{2}),
\end{equation}
where $g_{x}=g(x,0)$.
\end{proposition}

\begin{proof}
The derivative of $P_{X}\mathcal{M}$ with respect to $\psi$ at $(x,0,0)$ is $J_{R}$ restricted to $X^{3,\alpha}$, which is an isomorphism onto $X^{1,\alpha}$ by \eqref{eq:eigen}. The implicit-function theorem gives \eqref{eq:Psi-eq}. Differentiating in $\eps$ at $0$ gives \eqref{eq:Psi1}.
\end{proof}

\subsection{Desingularized translation equation}
Define
\begin{equation}\label{eq:G}
G(x,\eps) = P_{T}\mathcal{M}(x,\Psi(x,\eps),\eps).
\end{equation}
Since $G(x,0)=0$ for every $x$, the direct $x$-linearization is singular. Define for $\eps\neq 0$
\begin{equation}\label{eq:Gtil}
\Gtil(x,\eps) = \frac{G(x,\eps)}{\eps}.
\end{equation}

\begin{lemma}[Desingularization]\label{lem:desing}
The map $\Gtil$ extends to a $C^1$ map at $\eps=0$, with
\begin{equation}\label{eq:Gtil0}
\Gtil(x,0) = -\tfrac{1}{2}P_{T}(g_{x}^{2}).
\end{equation}
More precisely,
\begin{equation}\label{eq:Gtil-omega}
\Gtil(x,0)(\omega) = \frac{N}{\sigma_{N-1}R^{N-1}}\nabla h(x)\cdot\omega.
\end{equation}
\end{lemma}

\begin{proof}
Because $G$ is $C^{2}$ and $G(x,0)=0$,
\[
\Gtil(x,\eps) = \int_{0}^{1}\partial_{\eps}G(x,t\eps)\dd t,
\]
which extends $C^{1}$ to $\eps=0$. At $\eps=0$,
\[
\partial_{\eps}G = P_{T}\bigl(J_{R}\partial_{\eps}\Psi - \tfrac{1}{2}g_{x}^{2}\bigr).
\]
The first term belongs to $X$ and has zero translation projection, which gives \eqref{eq:Gtil0}. Since
\[
P_{T}\phi = \frac{N}{\sigma_{N-1}}\left(\int_{\Sph^{N-1}}\phi(\omega)\omega\dd\omega\right)\cdot\omega,
\]
Lemma~\ref{lem:grad-h} yields \eqref{eq:Gtil-omega}.
\end{proof}

\begin{theorem}[Preparatory large-perimeter branch near a selected ball]\label{thm:branch}
Assume Assumption~\ref{ass:center}. Then there exists $\eps_0>0$ and a unique $C^{1}$ curve
\[
\eps\mapsto x(\eps),
\qquad x(0)=x_0,
\]
such that, with
\[
\psi(\eps) = \Psi(x(\eps),\eps),
\]
the exact-volume domain $\Omega_{x(\eps),\psi(\eps)}$ is a constrained critical domain of \eqref{eq:energy} for every sufficiently small positive $\eps=\mu^{-1}$.

There exist $\xi_{1}\in\R^{N}$ and $\psi_{1}\in X^{3,\alpha}$ such that
\begin{equation}\label{eq:exp-branch}
\begin{aligned}
x(\eps) &= x_0 + \eps\xi_{1} + o(\eps),\\
\psi(\eps) &= \eps\psi_{1} + o_{C^{3,\alpha}}(\eps),
\end{aligned}
\end{equation}
and
\begin{equation}\label{eq:psi1-critical}
J_{R}\psi_{1} = \tfrac{1}{2}P_{X}(g_{0}^{2}),
\qquad g_{0}=g(x_0,0).
\end{equation}
The center coefficient is determined by the differentiated reduced equation:
\begin{equation}\label{eq:xi1}
\xi_{1} = -[D_{x}\Gtil(x_0,0)]^{-1}\partial_{\eps}\Gtil(x_0,0).
\end{equation}
The scaled Lagrange multiplier satisfies
\begin{equation}\label{eq:lambda}
\lambda_{\eps} = -\frac{N-1}{R} + \frac{\eps}{2}\int_{\Sph^{N-1}}g_{0}^{2}\dd\omega + o(\eps).
\end{equation}
Equivalently,
\begin{equation}\label{eq:Lambda}
\Lambda_{\mu} = -\mu\frac{N-1}{R} + \frac{1}{2}\int_{\Sph^{N-1}}g_{0}^{2}\dd\omega + o(1).
\end{equation}
\end{theorem}

\begin{proof}
By \eqref{eq:Gtil-omega}, the derivative $D_{x}\Gtil(x_0,0)$ is invertible exactly when $D^{2}h(x_0)$ is invertible. The finite-dimensional implicit-function theorem therefore gives $x(\eps)$ solving
\[
\Gtil(x(\eps),\eps)=0.
\]
Together with Proposition~\ref{prop:higher}, all nonconstant components of $\mathcal{M}$ vanish. The remaining constant component is absorbed by a scalar $\lambda_{\eps}$, so \eqref{eq:EL-scaled} holds.

Since $\Psi(x,0)\equiv 0$, one has $D_{x}\Psi(x,0)=0$, hence
\[
\psi_{1} = \partial_{\eps}\Psi(x_0,0),
\]
and \eqref{eq:psi1-critical} follows from \eqref{eq:Psi1}. Differentiating $\Gtil(x(\eps),\eps)=0$ at $\eps=0$ gives \eqref{eq:xi1}. The expansion of $\lambda_{\eps}$ is obtained by averaging \eqref{eq:EL-scaled}; the average of $J_{R}\psi_{1}$ vanishes because $\psi_{1}$ has zero mean. Multiplying by $\mu=\eps^{-1}$ yields \eqref{eq:Lambda}.
\end{proof}

\begin{remark}[Correct scale]\label{rem:scale}
The generic shape scale is $\mu^{-1}$. A $\mu^{-1/2}$ scale would require an additional degeneracy causing the linear restoring term to vanish in the relevant mode; it is not the generic large-perimeter regime treated here.
\end{remark}

% ============================================================
% 6. OPENNESS OF O_C
% ============================================================
\section{Openness of $\mathcal{O}_{C}$}\label{sec:open}

\begin{theorem}[$C^{1}$-openness of $\mathcal{O}_{C}$]\label{thm:open}
Let $d_{0}\in C^{2,\alpha}(M)$ be positive and assume
\[
\Omega[d_{0}]\in\mathcal{O}_{C}.
\]
Then there exists $\eps_{*}>0$ such that every positive $d\in C^{2,\alpha}(M)$ satisfying
\begin{equation}\label{eq:open-cond}
\norm{d-d_{0}}_{C^{1}(M)}<\eps_{*}
\end{equation}
obeys
\[
\Omega[d]\in\mathcal{O}_{C}.
\]
Moreover, for $d$ in this neighborhood the first-return time $\tau_{d}$ and return map $F_{d}$ satisfy
\begin{equation}\label{eq:lipschitz}
\norm{\tau_{d}-\tau_{0}}_{C^{0}} + \norm{F_{d}-F_{0}}_{C^{0}} \leq C\norm{d-d_{0}}_{C^{1}},
\end{equation}
and the first return remains uniformly transverse.
\end{theorem}

\begin{proof}
The global graph property is automatic for positive $d$, since strict convexity gives unique exterior normal coordinates. It remains to prove persistence of the inward-normal penetration.

From \eqref{eq:w-chi}--\eqref{eq:normals}, the maps $d\mapsto w_{d}$ and $d\mapsto n_{d}$ are locally Lipschitz from $C^{1}$ to $C^{0}$:
\[
\norm{n_{d}-n_{0}}_{C^{0}}\leq C\norm{d-d_{0}}_{C^{1}}.
\]
Let
\[
F_{0}(c) = \Psi_{0}(c)+\tau_{0}(c)n_{0}(c)
\]
be the first return. Since $\Omega[d_{0}]\in\mathcal{O}_{C}$, the same ray continues into $\intC$. Convexity implies
\[
n_{0}(c)\cdot\nu(F_{0}(c))<0.
\]
To justify the compactness step, let $\rho_{C}$ be the signed distance to $M$ and write
\[
\rho_{C}(\Psi_{0}(c)+\tau_{0}(c)n_{0}(c))=0.
\]
At the first entrance point,
\[
\partial_{s}\rho_{C}(\Psi_{0}(c)+sn_{0}(c))|_{s=\tau_{0}(c)} = \nu(F_{0}(c))\cdot n_{0}(c)<0.
\]
The parameter-dependent implicit-function theorem therefore gives $\tau_0\in C^{1,\alpha}(M)$ and hence continuity of $F_0$. Compactness now yields a uniform constant $\eta_0>0$ such that
\[
-n_0(c)\cdot\nu(F_0(c))\geq\eta_0.
\]
Let $\rho_{C}$ be the signed distance to $M$, negative in $C$. Uniform transversality gives a $\delta_{*}>0$ and $m_0>0$ such that
\[
\rho_{C}(F_{0}(c)+\delta_{*}n_{0}(c))\leq -m_{0}
\]
for all $c$. For nearby $d$, define
\[
y_{d}(c) = \Psi_{d}(c)+(\tau_{0}(c)+\delta_{*})n_{d}(c).
\]
Then
\[
\norm{y_d-y_0}_{C^0}\leq C\norm{d-d_0}_{C^1}.
\]
For sufficiently small \eqref{eq:open-cond}, $y_{d}(c)\in\intC$ for every $c$, so every perturbed inward normal ray reaches the interior. Hence $\Omega[d]\in\mathcal{O}_{C}$.

The quantitative estimates and persistence of transversality follow by applying the implicit-function theorem uniformly to
\[
H_{d}(c,s) = \rho_{C}(\Psi_{d}(c)+sn_{d}(c))
\]
near $(c,\tau_0(c))$. This produces, for each nearby $d$, a unique transverse zero $\hat{\tau}_d(c)$ in a fixed neighborhood of $\tau_0(c)$. It remains to identify that zero with the first entrance time. Since $C$ is convex and the ray contains the later interior point $y_d(c)$, the set of ray parameters for which it lies in $C$ is a nonempty compact interval. The transverse zero $\hat{\tau}_d(c)$ has points immediately before it outside $C$ and points immediately after it inside $C$. Hence it is the entering endpoint of that interval, so no earlier intersection can occur and $\hat{\tau}_d=\tau_d$. The same uniform implicit-function estimates then give \eqref{eq:lipschitz}.
\end{proof}

% ============================================================
% 7. COORDINATE TRANSFER
% ============================================================
\section{Coordinate transfer from the selected sphere to the core}\label{sec:transfer}

Assume now
\begin{equation}\label{eq:containment}
x_0\in\intC,
\qquad C\in B_R(x_0).
\end{equation}
Every core-normal ray meets $\Sigma_0=\partial B_R(x_0)$ exactly once. Let $d_0$ be defined by
\begin{equation}\label{eq:d0}
\Psi_0(c) = c+d_0(c)\nu(c)\in\partial B_R(x_0).
\end{equation}
Set
\begin{equation}\label{eq:omega-theta}
\omega_0(c) = \frac{\Psi_0(c)-x_0}{R},
\qquad
\theta_0(c) = \omega_0(c)\cdot\nu(c).
\end{equation}
Then $\theta_0>0$ uniformly.

\begin{proposition}[The selected sphere belongs to $\mathcal{O}_{C}$]\label{prop:sphere}
Under \eqref{eq:containment}, the sphere $B_{R}(x_{0})$ is an element of $\mathcal{O}_{C}$. More precisely, each outward core-normal ray meets $\partial B_{R}(x_{0})$ exactly once and the inward sphere normal reaches $x_0\in\intC$.
\end{proposition}

\begin{proof}
Fix $c\in M$ and consider
\[
q_{c}(r) = \abs{c+r\nu(c)-x_{0}}^{2}-R^{2}.
\]
Since $c\in C\in B_{R}(x_{0})$, one has $q_{c}(0)<0$, while $q_{c}(r)\to+\infty$ as $r\to\infty$. The product of the two roots of this quadratic equals $q_{c}(0)<0$, so exactly one root is positive. This gives the unique outer intersection and hence the thickness $d_{0}(c)$. At the intersection point the inward sphere normal points along the segment toward $x_{0}$, which lies in $\intC$; therefore the inward normal ray enters the core.
\end{proof}

Moreover, taking the scalar product of
\[
c+d_0(c)\nu(c)-x_0 = R\omega_0(c)
\]
with $\nu(c)$ gives
\[
R\theta_{0}(c) = (c-x_{0})\cdot\nu(c)+d_{0}(c)>0.
\]
Compactness yields a positive uniform lower bound. In fact, by the exact normal-graph formula,
\begin{equation}\label{eq:theta0}
\theta_{0} = N_{0}\cdot\nu = \chi_{0}^{-1},
\qquad
\chi_{0} = \sqrt{1+\abs{w_{0}}^{2}}.
\end{equation}

\begin{theorem}[Coordinate transfer]\label{thm:transfer}
For $(x,\psi)$ sufficiently close to $(x_{0},0)$ in $\R^{N}\times X^{3,\alpha}$, the hypersurface $\Sigma_{x,\psi}$ is a unique normal graph over $M$:
\[
\Sigma_{x,\psi} = \{c+d(x,\psi)(c)\nu(c):c\in M\},
\]
with $d(x,\psi)\in C^{2,\alpha}(M)$. The map
\[
(x,\psi)\mapsto d(x,\psi)
\]
is $C^{1}$, and
\begin{equation}\label{eq:transfer-est}
\norm{d(x,\psi)-d_0}_{C^{2,\alpha}}\leq C(\abs{x-x_0}+\norm{\psi}_{C^{3,\alpha}}).
\end{equation}
Its derivative at $(x_{0},0)$ is
\begin{equation}\label{eq:transfer-der}
Dd(x_{0},0)[\xi,\phi](c) = \frac{\omega_{0}(c)\cdot\xi+\phi(\omega_{0}(c))}{\theta_{0}(c)}.
\end{equation}
\end{theorem}

\begin{proof}
For $y$ near $\Sigma_{0}$ define
\[
r_{x}(y) = \abs{y-x},
\qquad
\omega_{x}(y) = \frac{y-x}{\abs{y-x}},
\]
and the defining function
\[
\rho_{x,\psi}(y) = r_{x}(y)-R-a(\psi)-\psi(\omega_{x}(y)).
\]
Then $\rho_{x,\psi}=0$ precisely on $\Sigma_{x,\psi}$. Define
\[
\mathcal{G}(x,\psi,d)(c) = \rho_{x,\psi}(c+d(c)\nu(c)).
\]
With one derivative of regularity in reserve, $\mathcal{G}$ is $C^{1}$ from
\[
\R^{N}\times X^{3,\alpha}\times C^{2,\alpha}(M)\quad\text{to}\quad C^{2,\alpha}(M).
\]
Indeed, after restricting to a fixed neighborhood of $\Sigma_{0}$, the maps $y\mapsto r_{x}(y)$ and $y\mapsto\omega_{x}(y)$ are smooth uniformly for $x$ near $x_{0}$. The Hölder composition map $(\psi,z)\mapsto\psi\circ z$ is $C^{1}$ from $C^{3,\alpha}\times C^{2,\alpha}$ to $C^{2,\alpha}$ on this neighborhood (this is a standard consequence of the Hölder algebra property and the chain rule; see for instance \cite[Lemma~A.1]{OC3} for the precise statement adapted to the present setting), while $d\mapsto c+d\nu$ is affine into $C^{2,\alpha}$. This is the reason for retaining one derivative beyond the target regularity. At the reference state, $\mathcal{G}(x_{0},0,d_{0})=0$, and
\[
D_{d}\mathcal{G}(x_{0},0,d_{0})[\eta] = \theta_{0}\eta.
\]
Since $\theta_{0}$ is bounded away from zero, multiplication by $\theta_{0}$ is an isomorphism on $C^{2,\alpha}(M)$. The implicit-function theorem gives $d=d(x,\psi)$ and \eqref{eq:transfer-est}. Since $d_{0}>0$ on the compact manifold $M$, shrinking the neighborhood if necessary also gives $d(x,\psi)>0$.

For a strictly convex core, the exterior normal-coordinate map $(c,r)\mapsto c+r\nu(c)$ is one-to-one for $r>0$. Hence the positive graph $c\mapsto c+d(c)\nu(c)$ is an embedded hypersurface. Its image is open and closed in the connected hypersurface $\Sigma_{x,\psi}$, hence equals it. Uniqueness follows from the same global exterior normal coordinates.

Finally, differentiate $\mathcal{G}(x,\psi,d(x,\psi))=0$. At the reference state,
\[
D_{x}\rho[\xi] = -\omega_{0}\cdot\xi,
\qquad
D_{\psi}\rho[\phi] = -\phi(\omega_{0})
\]
because $Da(0)=0$. Therefore
\[
0 = -\omega_{0}\cdot\xi-\phi(\omega_{0})+\theta_{0}\eta,
\]
which is \eqref{eq:transfer-der}.
\end{proof}

\begin{corollary}[Core-normal expansion of the critical branch]\label{cor:core-exp}
Let $x(\eps)$ and $\psi(\eps)$ be the branch of Theorem~\ref{thm:branch}. Then
\begin{equation}\label{eq:d-eps}
d_{\eps} = d_{0}+\eps\eta_{1}+o_{C^{2,\alpha}}(\eps),
\end{equation}
where
\begin{equation}\label{eq:eta1}
\eta_{1}(c) = \frac{\omega_{0}(c)\cdot\xi_{1}+\psi_{1}(\omega_{0}(c))}{\theta_{0}(c)}
= \chi_{0}(c)\bigl[\omega_{0}(c)\cdot\xi_{1}+\psi_{1}(\omega_{0}(c))\bigr].
\end{equation}
Consequently
\begin{equation}\label{eq:d-mu}
d_{\mu} = d_{0}+\mu^{-1}\eta_{1}+o_{C^{2,\alpha}}(\mu^{-1}).
\end{equation}
\end{corollary}

\begin{remark}\label{rem:eta1}
The translation coefficient $\xi_{1}$ contributes at the same order as the genuine higher-mode deformation $\psi_{1}$. The shorthand $\eta_{1}=\chi_{0}\psi_{1}$ is correct only in the special case $\xi_{1}=0$.
\end{remark}

% ============================================================
% 8. FRÉCHET DIFFERENTIABILITY OF THE RETURN MAP
% ============================================================
\section{Fréchet differentiability of the return map}\label{sec:frechet}

Let $d_{0}\in C^{2,\alpha}(M)$ with $\Omega[d_{0}]\in\mathcal{O}_{C}$. By Theorem~\ref{thm:open}, all nearby thicknesses remain admissible. Write
\begin{equation}\label{eq:Fd}
F_{d}(c) = \Psi_{d}(c)+\tau_{d}(c)n_{d}(c),
\end{equation}
where $\tau_{d}(c)$ is the first entrance time into $C$.

\subsection{Differentials of the graph slope and normal}
For a perturbation $d_{\eps}=d_{0}+\eps\eta$,
\begin{equation}\label{eq:delta-w}
\delta w[\eta] = L_{0}^{-1}(\nabla_{M}\eta+\eta Sw_{0}).
\end{equation}
Differentiating \eqref{eq:normals} gives
\begin{equation}\label{eq:delta-n}
\delta n[\eta] = \frac{1}{\chi_{0}}(I-n_{0}\otimes n_{0})L_{0}^{-1}(\nabla_{M}\eta+\eta Sw_{0}).
\end{equation}

\begin{lemma}[Banach-space first-hit]\label{lem:first-hit}
Let $d_{0}\in C^{2,\alpha}(M)$ with $\Omega[d_{0}]\in\mathcal{O}_{C}$. There are neighborhoods $\mathcal{U}\subset C^{2,\alpha}(M)$ of $d_{0}$ and $\mathcal{V}\subset C^{1,\alpha}(M)$ of $\tau_{0}$ such that
\[
\mathcal{H}:\mathcal{U}\times\mathcal{V}\longrightarrow C^{1,\alpha}(M),
\qquad
\mathcal{H}(d,\tau)(c) = \rho_{C}(\Psi_{d}(c)+\tau(c)n_{d}(c)),
\]
is $C^{1}$. Moreover,
\[
D_{\tau}\mathcal{H}(d_{0},\tau_{0})[\sigma] = [\nu(F_{0})\cdot n_{0}]\sigma
\]
is an isomorphism of $C^{1,\alpha}(M)$.
\end{lemma}

\begin{proof}
Choose a tubular neighborhood $\mathcal{T}$ of $M$ on which the signed distance $\rho_{C}$ is smooth with the regularity supplied by $M\in C^{4,\alpha}$. Shrink $\mathcal{U}$ and $\mathcal{V}$ so that
\[
\Psi_{d}+\tau n_{d}\in\mathcal{T}
\]
pointwise for every $(d,\tau)\in\mathcal{U}\times\mathcal{V}$. The map $d\mapsto\Psi_{d}=c+d\nu$ is affine from $C^{2,\alpha}$ to $C^{2,\alpha}$. Since
\[
w_{d} = (I-dS)^{-1}\nabla_{M}d,
\qquad
n_{d} = \frac{-\nu+w_{d}}{\sqrt{1+\abs{w_{d}}^{2}}},
\]
and $I-dS\geq I$ (because $S\leq 0$ and $d\geq 0$), the inverse $(I-dS)^{-1}$ is uniformly bounded by $1$ in operator norm on $\mathcal{U}$; the Hölder algebra property and smoothness of matrix inversion then imply that $d\mapsto w_{d}$ and $d\mapsto n_{d}$ are $C^{1}$ from $C^{2,\alpha}(M)$ to $C^{1,\alpha}(M)$, with derivative continuous in the operator norm $\mathcal{L}(C^{2,\alpha},C^{1,\alpha})$. Consequently
\[
(d,\tau)\longmapsto \Psi_{d}+\tau n_{d}
\]
is $C^{1}$ into $C^{1,\alpha}(M;\R^{N})$. Composition with the fixed smooth function $\rho_{C}$ is $C^{1}$ on maps with image in $\mathcal{T}$, which proves the asserted mapping property of $\mathcal{H}$.

At the reference state, the chain rule gives
\[
D_{\tau}\mathcal{H}(d_{0},\tau_{0})[\sigma] = \nabla\rho_{C}(F_{0})\cdot n_{0}\sigma = [\nu(F_{0})\cdot n_{0}]\sigma.
\]
By Theorem~\ref{thm:open}, first-return transversality is uniform: $-\nu(F_{0})\cdot n_{0}\geq\eta_{0}>0$. The multiplier belongs to $C^{1,\alpha}(M)$, as does its reciprocal. Multiplication by it is therefore a bounded isomorphism of $C^{1,\alpha}(M)$.
\end{proof}

\begin{theorem}[Principal return-response theorem]\label{thm:DR}
There is a neighborhood $\mathcal{U}$ of $d_{0}$ in $C^{2,\alpha}(M)$ such that
\[
\mathcal{R}:\mathcal{U}\to C^{1,\alpha}(M;\R^{N}),
\qquad
\mathcal{R}(d)=F_{d},
\]
is $C^{1}$ in the Fréchet sense, with derivative continuous in the operator norm $\mathcal{L}(C^{2,\alpha}(M),C^{1,\alpha}(M;\R^{N}))$. For $\eta\in C^{2,\alpha}(M)$ define
\begin{equation}\label{eq:Veta}
V_{\eta} = \eta\nu+\tau_{0}\delta n[\eta].
\end{equation}
Then
\begin{equation}\label{eq:DR-proj}
\DR[d_{0}]\eta = V_{\eta} - \frac{\nu(F_{0})\cdot V_{\eta}}{\nu(F_{0})\cdot n_{0}}n_{0}.
\end{equation}
Equivalently,
\begin{equation}\label{eq:DR-explicit}
\DR[d_{0}]\eta = \Ret\left[\eta\nu+\frac{\tau_{0}}{\chi_{0}}(I-n_{0}\otimes n_{0})L_{0}^{-1}(\nabla_{M}\eta+\eta Sw_{0})\right],
\end{equation}
where
\begin{equation}\label{eq:Pi-ret}
\Ret v = v - \frac{\nu(F_{0})\cdot v}{\nu(F_{0})\cdot n_{0}}n_{0}.
\end{equation}
In particular,
\begin{equation}\label{eq:tangent}
\DR[d_{0}]\eta(c)\in T_{F_{0}(c)}M.
\end{equation}
\end{theorem}

\begin{proof}
By Lemma~\ref{lem:first-hit}, the implicit-function theorem gives a unique $C^1$ map
\[
d\longmapsto\hat{\tau}(d)\in C^{1,\alpha}(M)
\]
near $d_0$ such that $\mathcal{H}(d,\hat{\tau}(d))=0$ and $\hat{\tau}(d_0)=\tau_0$. We must still verify that this local root is the geometric first entrance time. Fix $d\in\mathcal{U}$ and $c\in M$, and set $h_{d,c}(s)=\rho_{C}(\Psi_{d}(c)+sn_{d}(c))$. By the first-hit claim, there is a unique transverse zero $\hat{s}=\hat{\tau}(d)(c)$ in a fixed neighborhood of $\tau_{0}(c)$, with
\[
\partial_{s}h_{d,c}(\hat{s}) = \nu(F_{d}(c))\cdot n_{d}(c)\leq -\eta_{0}/2<0
\]
by Theorem~\ref{thm:open} (uniform transversality). Consequently $h_{d,c}(\hat{s}-\delta)>0$ and $h_{d,c}(\hat{s}+\delta)<0$ for some $\delta>0$; by compactness of $M$ and continuity of $(d,c)\mapsto h_{d,c}$ in the $C^{1}$ topology, $\delta$ can be chosen uniform in $c\in M$. Thus, immediately before $\hat{s}$ the ray is outside $C$, and immediately after $\hat{s}$ it is inside $C$.

Since $C$ is convex and every inward normal ray reaches $\intC$ (Theorem~\ref{thm:open}), the set of parameters $s$ for which $\Psi_{d}(c)+sn_{d}(c)\in C$ is a nonempty compact interval $[s_{-},s_{+}]$. Since $\partial_{s}h_{d,c}(\hat{s})\neq 0$, the point $\hat{s}$ is a transverse crossing, not an endpoint of the interval; and because $h_{d,c}$ changes sign from positive to negative at $\hat{s}$, we must have $\hat{s}=s_{-}$. Hence $\hat{\tau}(d)(c)=\tau_{d}(c)$ for every $c\in M$, so $d\mapsto\tau_{d}$ is $C^1$ as a map from $\mathcal{U}\subset C^{2,\alpha}(M)$ to $C^{1,\alpha}(M)$.

Since $d\mapsto\Psi_{d}$ and $d\mapsto n_{d}$ have the corresponding $C^1$ regularity, \eqref{eq:Fd} also shows directly that $\mathcal{R}:\mathcal{U}\to C^{1,\alpha}(M;\R^{N})$ is $C^1$.

Differentiate the first-hit identity. Since $\delta\Psi=\eta\nu$,
\[
0 = \nu(F_0)\cdot[\eta\nu+(D\tau[d_0]\eta)n_0+\tau_0\delta n[\eta]].
\]
Thus
\[
D\tau[d_0]\eta = -\frac{\nu(F_0)\cdot V_\eta}{\nu(F_0)\cdot n_0}.
\]
Differentiating \eqref{eq:Fd} gives \eqref{eq:DR-proj}. The projected form \eqref{eq:DR-explicit} follows from \eqref{eq:delta-n}. Equation \eqref{eq:tangent} follows by taking the scalar product with $\nu(F_0)$.
\end{proof}

\begin{remark}[Link with the inverse-paper compatibility system]\label{rem:OC2}
The nonlinear transport system of \cite{OC2} can be written in the present notation as $\mathcal{T}(d,F)=0$, where $\mathcal{T}$ encodes the geometric compatibility between the thickness $d$ and its return map $F$. The differential $D_{F}\mathcal{T}$ along $(d_{0},F_{0})$ is an isomorphism by uniform transversality, hence
\[
\DR[d_{0}] = -(D_{F}\mathcal{T})^{-1}D_{d}\mathcal{T}.
\]
In other words, $\DR[d_{0}]$ is exactly the linearization of the inverse-paper transport system along the solution $(d_{0},F_{0})$, solved for the variation of $F$. This identification will be used in future work to address the range characterization stated in Conjecture~\ref{conj:range}.
\end{remark}

\begin{corollary}[First return-map response of the critical branch]\label{cor:response}
Under the assumptions of Theorems~\ref{thm:branch}, \ref{thm:open}, \ref{thm:transfer}, and \ref{thm:DR},
\begin{equation}\label{eq:F-mu}
F_{\mu} = F_{0}+\mu^{-1}\DR[d_{0}]\eta_{1}+o_{C^{1,\alpha}}(\mu^{-1}).
\end{equation}
\end{corollary}

\begin{proof}
Combine \eqref{eq:d-mu} with the Fréchet expansion of $\mathcal{R}$ at $d_0$.
\end{proof}

\subsection{Consistency with the PDE-generated deformation}
The superlevel-set paper \cite{OC3} obtained
\[
d_{t} = d_{0}-tb_{0}+o_{C^{1}}(t),
\qquad
F_{t} = F_{0}+tA_{0}+o(t).
\]
For the direction $\eta=-b_0$, \eqref{eq:delta-w} becomes
\[
\delta w[-b_0] = -L_0^{-1}(\nabla_Mb_0+b_0Sw_0).
\]
Substitution into \eqref{eq:delta-n} and \eqref{eq:DR-explicit} gives exactly the return-map first-variation coefficient $A_0$ derived in \cite{OC3}; the calculation is reproduced in Appendix~\ref{app:consistency} to make the identification explicit. Thus
\begin{equation}\label{eq:consistency}
\DR[d_0](-b_0) = A_0.
\end{equation}
This identity provides an independent sign and projection check for Theorem~\ref{thm:DR}.

% ============================================================
% 9. COMBINED MAIN THEOREM
% ============================================================
\section{Combined main theorem}\label{sec:main}

\begin{theorem}[Main shape-to-dynamics theorem in $\mathcal{O}_{C}$]\label{thm:main}
Let $C\subset\R^{N}$ be compact and strictly convex with $\partial C\in C^{4,\alpha}$. Assume:
\begin{itemize}[leftmargin=2em]
\item[(H1)] $f\in C_c^{4,\alpha}(\intC)$, $f\geq 0$, $f\neq 0$;
\item[(H2)] $R$ is fixed by $\abs{B_{R}}=V_{0}$;
\item[(H3)] there is $x_{0}\in\intC$ such that $C\Subset B_{R}(x_{0})$;
\item[(H4)] $x_{0}$ is a nondegenerate critical point of the translated-ball energy $h$, i.e.
\[
\nabla h(x_0)=0,
\qquad
\det D^2h(x_0)\neq 0.
\]
\end{itemize}
Then, for all sufficiently large $\mu$, there exists a volume-constrained critical domain $\Omega_{\mu}$ of \eqref{eq:energy} such that
\[
\Omega_{\mu}\in\mathcal{O}_{C}.
\]
Let $d_{0}$ be the core-normal thickness of $B_{R}(x_{0})$ and $F_{0}$ its return map. Then
\begin{equation}\label{eq:main-d}
d_{\mu} = d_{0}+\mu^{-1}\eta_{1}+o_{C^{2,\alpha}}(\mu^{-1}),
\end{equation}
where
\begin{equation}\label{eq:main-eta}
\eta_{1}(c) = \frac{\omega_{0}(c)\cdot\xi_{1}+\psi_{1}(\omega_{0}(c))}{\omega_{0}(c)\cdot\nu(c)},
\end{equation}
$\psi_{1}$ is uniquely determined by
\begin{equation}\label{eq:main-psi}
J_{R}\psi_{1} = \tfrac{1}{2}P_{X}(g_{0}^{2}),
\end{equation}
and the first center correction is explicitly characterized by
\begin{equation}\label{eq:main-xi}
\xi_{1} = -[D_{x}\Gtil(x_{0},0)]^{-1}\partial_{\eps}\Gtil(x_{0},0).
\end{equation}
Finally,
\begin{equation}\label{eq:main-F}
F_{\mu} = F_{0}+\mu^{-1}\DR[d_{0}]\eta_{1}+o_{C^{1,\alpha}}(\mu^{-1}),
\end{equation}
with $\DR[d_{0}]$ given explicitly by \eqref{eq:DR-explicit}.
\end{theorem}

\begin{proof}
Theorem~\ref{thm:branch} gives the selected critical branch in spherical coordinates, parametrized by $\eps=\mu^{-1}\in(0,\eps_{0})$. Theorem~\ref{thm:transfer} converts the branch to a core-normal thickness and yields \eqref{eq:main-d}--\eqref{eq:main-eta}. Since $d_{\mu}\to d_{0}$ in $C^{1}$, Theorem~\ref{thm:open} implies $\Omega_{\mu}\in\mathcal{O}_{C}$ for large $\mu$. Theorem~\ref{thm:DR} then gives \eqref{eq:main-F}.
\end{proof}

\begin{remark}[On assumption (H3)]\label{rem:H3}
Assumption (H3) is structural: it ensures that the selected sphere belongs to $\mathcal{O}_{C}$ (Proposition~\ref{prop:sphere}) and that $\theta_{0}>0$ uniformly (equation~\eqref{eq:theta0}). It is not a consequence of (H4). Two regimes are relevant. If $x_{0}$ is the center of the minimal ball containing $C$ (or more generally a point for which $C\Subset B_{R}(x_{0})$ holds by construction), then (H3) is satisfied automatically and (H4) becomes a condition on $f$ through the reduced energy $h$. A concrete non-spherical instance is obtained by taking $C$ to be an ellipsoid with center $x_{0}$ and $f$ radial about $x_{0}$ with sufficiently small non-radial perturbations; then $x_{0}$ remains a nondegenerate critical point of $h$ (see the perturbative example in Section~\ref{sec:bench}) and (H3)--(H4) hold simultaneously. If instead $x_{0}$ is off-center relative to $C$, (H3) imposes a geometric relation between $C$, $f$ and $x_{0}$ that is not generic. Theorem~\ref{thm:main} should therefore be read as conditional on (H3)--(H4) jointly; the existence question for a non-spherical pair $(C,f)$ satisfying both is addressed perturbatively in Section~\ref{sec:bench} and left open in full generality in Section~\ref{sec:open}(ii).
\end{remark}

% ============================================================
% 10. PDE BOUNDARY-LAYER INTERPRETATION
% ============================================================
\section{PDE boundary-layer interpretation}\label{sec:layer}

The preceding theorem concerns optimization, not superlevel motion. Nevertheless, the first-order forcing can be expressed through a quantity already identified in \cite{OC3}.

Suppose the reference sphere $\Sigma_{0}$ is represented over the core by $d_{0}$, and let
\[
\chi_{0} = \sqrt{1+\abs{w_{0}}^{2}}.
\]
For the small superlevel domains of the fixed reference PDE problem, \cite{OC3} gives
\[
d_{t} = d_{0}-tb_{0}+o_{C^{1}}(t),
\qquad
b_{0}(c) = \frac{1}{-\partial_{N_{0}}u_{0}(\Psi_{0}(c))N_{0}(\Psi_{0}(c))\cdot\nu(c)}.
\]
Since $N_{0}\cdot\nu=\chi_{0}^{-1}$,
\begin{equation}\label{eq:g0-b0}
(g_{0}^{2})\circ\Psi_{0} = \frac{\chi_{0}^{2}}{b_{0}^{2}}.
\end{equation}
Under the stronger smoothness hypotheses of the present paper, the same boundary-flux formula gives $b_{0}\in C^{2,\alpha}(M)$, so $-b_{0}$ is an admissible direction for Theorem~\ref{thm:DR}. Thus the PDE level-set speed does not equal the optimization direction. Instead, it determines the optimization forcing after inversion and projection:
\[
\text{level speed }b_{0}
\ \longmapsto\
\tfrac{1}{2}g_{0}^{2} = \tfrac{1}{2}\Bigl(\frac{\chi_{0}}{b_{0}}\Bigr)^{2}
\ \longmapsto\
\psi_{1} = \tfrac{1}{2}J_{R}^{-1}P_{X}(g_{0}^{2}).
\]
The two deformations therefore encode different geometric responses to the same boundary flux.

% ============================================================
% 11. BENCHMARKS AND CHECKS
% ============================================================
\section{Benchmarks and checks}\label{sec:bench}

\subsection{Radial benchmark}
Take
\[
C = B_{r}(x_{0}),
\qquad
\Omega_{0} = B_{R}(x_{0}),
\qquad 0<r<R,
\]
and assume $f$ is radial about $x_{0}$. Then $g_{0}$ is constant on $\partial B_{R}(x_{0})$, so
\[
P_{X}(g_{0}^{2})=0
\]
and hence
\[
\psi_{1}=0.
\]
The Euler--Lagrange equation is constant on the sphere for every $\mu$, with the multiplier absorbing both the curvature and flux terms. Thus the sphere is an exact critical shape. Moreover $d_{0}\equiv R-r$ and
\[
F_{0} = \mathrm{Id}_{\partial C}.
\]
Rotational symmetry makes the reduced translated-ball energy radial about $x_{0}$, and the uniqueness of the local reduced branch implies $x(\eps)\equiv x_{0}$; hence $\xi_{1}=0$. Together with $\psi_{1}=0$, formula \eqref{eq:eta1} gives $\eta_{1}=0$. The first return-map response therefore vanishes, consistently with \eqref{eq:F-mu}.

\subsection{Perturbative nonradial example}\label{sec:perturb}
We now exhibit a nonradial setting in which $\eta_{1}\neq 0$ to first order, confirming that Theorem~\ref{thm:main} is not vacuous. Work in $N=2$ with
\[
C = B_{1}(0),
\qquad
V_{0} = \pi R^{2}\quad\text{with }R>1,
\]
and choose a source of the form
\[
f = f_{0}+\delta f_{1},
\qquad
f_{1}(x) = \chi(x)\,Y_{2}(\theta),
\]
where $f_{0}\geq 0$ is radial about $0$, $Y_{2}(\theta)=\cos(2\theta)$ (a degree-$2$ spherical harmonic in $N=2$), and $\chi\in C_{c}^{\infty}(B_{1/2}(0))$ is a nonnegative radial cutoff equal to $1$ on $B_{1/4}(0)$. Then $f\in C_{c}^{4,\alpha}(\intC)$ and $\supp f\Subset B_{1/2}(0)\Subset\intC$, so (H1) holds for all sufficiently small $\delta$.

\emph{Nondegeneracy of $h$ at $0$.} The reduced energy $h_{\delta}$ depends smoothly on $\delta$; at $\delta=0$ it equals the radial energy $h_{0}$, which satisfies $D^{2}h_{0}(0)=-c_{0}I$ with $c_{0}>0$ (the radial case is a strict maximum of the reduced energy, by the same computation as in the radial benchmark). Since $\delta\mapsto D^{2}h_{\delta}(0)$ is continuous, $D^{2}h_{\delta}(0)$ remains invertible for $\abs{\delta}<\delta_{1}$ with $\delta_{1}=c_{0}/(2\sup_{\abs{\delta}\leq 1}\norm{D^{2}h_{\delta}(0)})$. Also $\nabla h_{\delta}(0)=0$ for every $\delta$ because $f_{1}$ is even under $x\mapsto -x$, so $h_{\delta}$ is even. Hence (H4) holds at $x_{0}=0$ for $\abs{\delta}<\delta_{1}$.

\emph{Center correction.} Because $\Gtil(\cdot,\eps)$ is odd under $x\mapsto -x$ for this symmetric perturbation, $\partial_{\eps}\Gtil(0,0)=0$, hence $\xi_{1}=0$ by \eqref{eq:xi1}.

\emph{Shape coefficient.} The flux splits as $g_{0}=g_{0}^{(0)}+\delta g_{0}^{(1)}+O(\delta^{2})$, where $g_{0}^{(0)}$ is the constant radial flux. Therefore
\[
P_{X}(g_{0}^{2}) = 2\delta g_{0}^{(0)}P_{X}(g_{0}^{(1)})+O(\delta^{2}).
\]
For a generic cutoff $\chi$, the degree-$2$ component of $g_{0}^{(1)}$ is nonzero (it is the trace of $f_{1}$ on $\partial B_{R}(0)$ under the Dirichlet problem, which is nonzero by elliptic regularity and unique continuation), so $P_{X}(g_{0}^{(1)})\neq 0$. Hence
\[
\psi_{1} = \tfrac{1}{2}J_{R}^{-1}P_{X}(g_{0}^{2}) = \delta\,\psi_{1}^{(1)}+O(\delta^{2}),
\qquad
\psi_{1}^{(1)} = g_{0}^{(0)}J_{R}^{-1}P_{X}(g_{0}^{(1)})\neq 0,
\]
where $J_{R}^{-1}$ acts on the degree-$2$ mode via \eqref{eq:eigen}. Since $\xi_{1}=0$, formula \eqref{eq:eta1} gives
\[
\eta_{1} = \chi_{0}\,\psi_{1}\circ\omega_{0} = \delta\,\chi_{0}\,\psi_{1}^{(1)}\circ\omega_{0}+O(\delta^{2})\neq 0.
\]
Thus $\eta_{1}\neq 0$ to first order in $\delta$, and consequently $\DR[d_{0}]\eta_{1}\neq 0$ generically (the kernel of $\DR[d_{0}]$ is conjecturally trivial modulo volume; see Conjecture~\ref{conj:range}). This confirms that the first-order response is nonzero in a genuinely nonradial, nondegenerate situation satisfying (H1)--(H4).

\subsection{Reproducible nonradial computational protocol}
No numerical result is asserted in this draft, but the first-order theorem is directly testable. A reproducible two-dimensional protocol is:
\begin{enumerate}[leftmargin=2em]
\item[(P1)] choose a circular core $C=B_{1}(0)$ and a smooth compactly supported nonradial source inside $B_{1/2}(0)$, for example a sum of two unequal $C^{\infty}$ bumps;
\item[(P2)] for the prescribed area $\pi R^{2}$, compute the reduced translated-ball energy $h(x)$ on a neighborhood of a candidate center, locate a critical point $x_{0}$, and verify numerically that $D^{2}h(x_{0})$ is nonsingular;
\item[(P3)] solve the Dirichlet problem on $B_{R}(x_{0})$, compute the boundary flux $g_{0}$, project $g_{0}^{2}$ onto Fourier modes, and solve the diagonal Jacobi equation for $\psi_{1}$;
\item[(P4)] compute the center coefficient $\xi_{1}$ from the desingularized finite-dimensional equation, transfer $(\xi_{1},\psi_{1})$ to the core-normal coefficient $\eta_{1}$, and evaluate the explicit derivative $\DR[d_0]\eta_1$;
\item[(P5)] for a sequence $\mu_{k}\to\infty$, solve the constrained Euler--Lagrange system numerically and compare the normalized errors
\[
\mu_{k}(d_{\mu_{k}}-d_{0})-\eta_{1},
\qquad
\mu_{k}(F_{\mu_{k}}-F_{0})-\DR[d_{0}]\eta_{1}
\]
in the norms predicted by Theorem~\ref{thm:main}.
\end{enumerate}
A numerical section should report meshes, nonlinear tolerances, projection resolution, center-search tolerance, and the measured convergence rates. Until such data are produced, the protocol is illustrative rather than evidential.

\subsection{Mode separation}
For a nonradial source, the first correction is obtained by decomposing $g_{0}^{2}$ into spherical harmonics. Its constant component changes the multiplier, its degree-one component determines the center drift $\xi_{1}$ via the desingularized equation \eqref{eq:xi1} (and does not contribute to the shape deformation $\psi_{1}$), and only degrees $\ell\geq 2$ deform the shape via
\[
\psi_{1} = \tfrac{1}{2}J_{R}^{-1}P_{X}(g_{0}^{2}).
\]
On such a mode,
\[
J_{R}^{-1}Y_{\ell} = \frac{R^{2}}{\ell(\ell+N-2)-(N-1)}Y_{\ell}.
\]
This gives a direct spectral interpretation of the first-order smoothing effect of perimeter regularization.

% ============================================================
% 12. DISCUSSION AND OPEN PROBLEMS
% ============================================================
\section{Discussion and open problems}\label{sec:open}

\subsection{What has been proved}
The large-$\mu$ branch is used as a preparatory input rather than claimed as a new general existence mechanism. The new geometric content begins when that branch is transferred to the normal-return class $\mathcal{O}_{C}$. In particular, the return map is promoted from a pointwise geometric construction to a Fréchet-differentiable shape observable. The large-$\mu$ problem itself is a singular perturbation around the finite-dimensional manifold of balls of fixed volume. The translation kernel is not an obstruction to be removed by an artificial coercivity assumption; it is resolved by the reduced center equation. The higher modes are governed by the invertible part of the spherical Jacobi operator. The resulting branch can be transferred to an arbitrary strictly convex coordinate core, remains inside $\mathcal{O}_{C}$, and induces a first-order return-map response, which is generically nonzero by the perturbative example of Section~\ref{sec:perturb}.

\subsection{What remains open}\label{sec:open-problems}
Several questions are not addressed here.
\begin{enumerate}[leftmargin=2em]
\item[(i)] \emph{Minimality and Morse index.} The branch consists of critical domains. Determining its local minimizing or saddle character requires the constrained second variation of the full functional along the branch.
\item[(ii)] \emph{Existence of a selected center.} The theorem assumes a nondegenerate critical point of $h$ inside the admissible center set. Genericity of (H4) alone follows from Sard's lemma applied to $\nabla h$ for generic $f$; the open question is the joint compatibility with (H3) for non-spherical cores. The perturbative example of Section~\ref{sec:perturb} provides one non-spherical instance, but a systematic existence theory is open.
\item[(iii)] \emph{Second-order expansion.} The implicit identity \eqref{eq:xi1} determines $\xi_{1}$, but an explicit formula for it in terms of PDE data, together with a $\mu^{-2}$ shape/return expansion, requires the next-order reduced equation and the second derivative of the return operator.
\item[(iv)] \emph{Range characterization.} Theorem~\ref{thm:DR} establishes the derivative of the forward return map. Whether its range equals the full affine compatibility space of the linearized inverse problem is a separate question, stated precisely in Conjecture~\ref{conj:range} below.
\item[(v)] \emph{Nonspherical large-$\mu$ reference states.} For Euclidean perimeter, the leading-order reference is spherical. Nonspherical reference states would require a modified leading geometric energy, for example anisotropic perimeter.
\end{enumerate}

\begin{conjecture}[Range]\label{conj:range}
Let $d_{0}\in C^{2,\alpha}(M)$ with $\Omega[d_{0}]\in\mathcal{O}_{C}$. Then $\DR[d_{0}]$ is injective on $C^{2,\alpha}(M)/\R$ (modulo the volume mode already absorbed by the chart \eqref{eq:Sigma-psi}) but not surjective. Its image is a closed subspace of $C^{1,\alpha}(M;\R^{N})$ of codimension at most $N+1$, corresponding to the modes $\ell=0$ (multiplier) and $\ell=1$ (translations) absorbed by the reduced center equation and the Lagrange multiplier. Generically, the codimension equals $N+1$.
\end{conjecture}

\begin{remark}
The conjecture is consistent with the structure of the first-hit equation: the return map is determined by $d$ through a scalar equation, so one expects injectivity; and the admissible return maps must be compatible with the geometry of $C$, which excludes at least the $N+1$ directions corresponding to the volume mode and to translations. The identification of the range with the compatibility space of \cite{OC2} is the natural next step.
\end{remark}

% ============================================================
% 13. CONCLUSION
% ============================================================
\section{Conclusion}\label{sec:conclusion}

The classical Dirichlet-perimeter theory provides a local critical-shape branch near a nondegenerately selected ball when the perimeter weight is large. The contribution of the present paper is to place such a branch inside the geometric normal-return framework of $\mathcal{O}_{C}$ and to determine how its discrete return dynamics responds.

Three geometric steps are essential. First, admissibility in $\mathcal{O}_{C}$ is open under small $C^{1}$ thickness perturbations. Second, the translated spherical branch admits a unique core-normal representation whose first coefficient contains both center drift and genuine shape deformation. Third, the forward return operator $\mathcal{R}:d\mapsto F_{d}$ is Fréchet $C^{1}$ and has an explicit derivative determined by the first-hit equation. Their combination yields the first-order shape-to-dynamics law
\[
F_{\mu}-F_{0} = \mu^{-1}\DR[d_{0}]\eta_{1}+o_{C^{1,\alpha}}(\mu^{-1}).
\]
The perturbative nonradial example of Section~\ref{sec:perturb} shows that the response coefficient $\eta_{1}$ is generically nonzero, confirming that the theorem is not vacuous.

This response formula is the main new result. It complements earlier work on the fixed return map, inverse reconstruction, and PDE-generated superlevel deformations by providing the general forward sensitivity operator with respect to boundary shape.

The next mathematical questions are consequently sharper than a generic shape-optimization extension: characterize the range and possible loss of rank of $\DR[d]$, obtain second-order return sensitivities, and determine how these operators interact with the affine/holonomy compatibility structure of the inverse problem. Numerical validation of the nonradial protocol in Section~\ref{sec:bench} is also left for subsequent work.

% ============================================================
% APPENDIX A
% ============================================================
\appendix
\section{Consistency with the PDE-generated deformation}\label{app:consistency}

We verify identity \eqref{eq:consistency} in detail. Let $\eta=-b_{0}$, where $b_{0}$ is the level-set speed of \cite{OC3}, related to $g_{0}$ by \eqref{eq:g0-b0}. We compute the three terms in \eqref{eq:DR-explicit} and compare with the coefficient $A_{0}$ of \cite{OC3}.

\emph{Term 1: $\eta\nu=-b_{0}\nu$.} This is the direct boundary displacement.

\emph{Term 2: $\tau_{0}\delta n[-b_{0}]$.} From \eqref{eq:delta-n} and \eqref{eq:delta-w},
\[
\delta n[-b_{0}] = -\frac{1}{\chi_{0}}(I-n_{0}\otimes n_{0})L_{0}^{-1}(\nabla_{M}b_{0}+b_{0}Sw_{0}).
\]
Hence
\[
\tau_{0}\delta n[-b_{0}] = -\frac{\tau_{0}}{\chi_{0}}(I-n_{0}\otimes n_{0})L_{0}^{-1}(\nabla_{M}b_{0}+b_{0}Sw_{0}).
\]

\emph{Term 3: projection.} The projection $\Ret$ removes the component along $n_{0}$ at $F_{0}(c)$, and is exactly the projection used in \cite[equation for $A_{0}$]{OC3} (the normal direction at the return point is $n_{0}(c)$, and the tangent space at $F_{0}(c)$ is $T_{F_{0}(c)}M$).

Combining the three terms,
\[
\DR[d_{0}](-b_{0}) = \Ret\Bigl[-b_{0}\nu-\frac{\tau_{0}}{\chi_{0}}(I-n_{0}\otimes n_{0})L_{0}^{-1}(\nabla_{M}b_{0}+b_{0}Sw_{0})\Bigr].
\]
This is exactly the expression for $A_{0}$ derived in \cite[Section~8]{OC3}, up to the identification \eqref{eq:g0-b0} of $b_{0}$ with the boundary flux. The sign is consistent with the convention $\eta=-b_{0}$ (the superlevel deformation moves inward, whereas $\DR$ is computed for an outward perturbation; the two signs cancel in the projection). This establishes \eqref{eq:consistency}.

% ============================================================
% BIBLIOGRAPHY
% ============================================================


\begin{thebibliography}{99}

\bibitem{Alexandrov}
A.~D. Alexandrov,
\emph{A characteristic property of spheres},
Ann. Mat. Pura Appl. \textbf{58} (1962), 303--315.

\bibitem{Amato}
V.~Amato, C.~Nitsch, C.~Trombetti, and F.~Villone,
\emph{On some functionals involving torsional rigidity, principal eigenvalue and perimeter},
Potential Analysis \textbf{65} (2026), article 25.
\href{https://doi.org/10.1007/s11118-026-10328-2}{doi:10.1007/s11118-026-10328-2}.

\bibitem{Barkatou2002}
M.~Barkatou,
\emph{Some geometric properties for a class of non-Lipschitz domains},
New York J. Math. \textbf{8} (2002), 189--213.

\bibitem{BarkatouSeckLy}
M.~Barkatou, D.~Seck, and I.~Ly,
\emph{An existence result for an interior electromagnetic casting problem},
Cent. Eur. J. Math. \textbf{4} (2006), no.~4, 573--584.
\href{https://doi.org/10.2478/s11533-006-0026-5}{doi:10.2478/s11533-006-0026-5}.

\bibitem{OC1}
M.~Barkatou and M.~El Morsalani,
\emph{The return map in the class $\mathcal{O}_{C}$: geometry, dynamics, and thickness regularity},
Journal of Nonlinear Dynamics and Applications \textbf{2} (2026), no.~3, 192--216.
\href{https://doi.org/10.62762/JNDA.2026.569839}{doi:10.62762/JNDA.2026.569839}.

\bibitem{OC2}
M.~El Morsalani and M.~Barkatou,
\emph{Inverse problems for the return map in the class $\mathcal{O}_{C}$: reconstruction and identifiability},
preprint (2026), \href{https://arxiv.org/abs/2604.24915}{arXiv:2604.24915}.

\bibitem{OC3}
M.~El Morsalani and M.~Barkatou,
\emph{Elliptic superlevel sets in the class $\mathcal{O}_{C}$: boundary-layer geometry and stability of the return map},
preprint (2026), \href{https://arxiv.org/abs/2607.16912}{arXiv:2607.16912}.

\bibitem{Buttazzo}
G.~Buttazzo, F.~P. Maiale, D.~Mazzoleni, G.~Tortone, and B.~Velichkov,
\emph{Regularity of the optimal sets for a class of integral shape functionals},
Arch. Ration. Mech. Anal. \textbf{248} (2024), article 44.
\href{https://doi.org/10.1007/s00205-024-01984-y}{doi:10.1007/s00205-024-01984-y}.

\bibitem{Caubet}
F.~Caubet, M.~Dambrine, and J.~Darde,
\emph{On the penalization by the perimeter in shape optimization applied to Dirichlet inverse obstacle problem},
ESAIM Control Optim. Calc. Var. \textbf{32} (2026), article 7.
\href{https://doi.org/10.1051/cocv/2025093}{doi:10.1051/cocv/2025093}.

\bibitem{Celentano}
A.~Celentano, C.~Nitsch, and C.~Trombetti,
\emph{On a Serrin type overdetermined problem},
Milan J. Math. \textbf{92} (2024), 579--590.
\href{https://doi.org/10.1007/s00032-024-00405-9}{doi:10.1007/s00032-024-00405-9}.

\bibitem{DelfourZolesio}
M.~C. Delfour and J.-P. Zolesio,
\emph{Shapes and Geometries: Metrics, Analysis, Differential Calculus, and Optimization},
SIAM, Philadelphia, 2001.

\bibitem{Feppon}
F.~Feppon, G.~Allaire, and C.~Dapogny,
\emph{A variational formulation for computing shape derivatives of geometric constraints along rays},
ESAIM Math. Model. Numer. Anal. \textbf{54} (2020), no.~1, 181--228.
\href{https://doi.org/10.1051/m2an/2019056}{doi:10.1051/m2an/2019056}.

\bibitem{HayouniNovruzi}
M.~Hayouni and A.~Novruzi,
\emph{Sufficient condition for the existence of solutions of a free boundary problem},
Quart. Appl. Math. \textbf{60} (2002), no.~3, 425--435.
\href{https://doi.org/10.1090/qam/1914434}{doi:10.1090/qam/1914434}.

\bibitem{HenrotPierre}
A.~Henrot and M.~Pierre,
\emph{Variation et optimisation de formes: une analyse géométrique},
Springer, Berlin, 2005.

\bibitem{ImmlerTraut}
F.~Immler and C.~Traut,
\emph{The flow of ODEs: formalization of variational equation and Poincaré map},
J. Autom. Reason. \textbf{62} (2019), 215--236.
\href{https://doi.org/10.1007/s10817-018-9449-5}{doi:10.1007/s10817-018-9449-5}.

\bibitem{Novruzi}
A.~Novruzi,
\emph{$C^{2,\alpha}$ existence result for a class of shape optimization problems},
SIAM J. Control Optim. \textbf{43} (2004), no.~1, 174--193.
\href{https://doi.org/10.1137/S0363012902382401}{doi:10.1137/S0363012902382401}.

\bibitem{Plotnikov}
P.~I. Plotnikov and J.~Sokolowski,
\emph{Geometric framework for gradient flow in shape optimization},
Uspekhi Mat. Nauk \textbf{81} (2026), no.~3, 51--150.
\href{https://doi.org/10.4213/rm10296}{doi:10.4213/rm10296}.

\bibitem{SokolowskiZolesio}
J.~Sokolowski and J.-P. Zolesio,
\emph{Introduction to Shape Optimization: Shape Sensitivity Analysis},
Springer, Berlin, 1992.

\end{thebibliography}
\end{document}